\documentclass[12pt]{article}
\usepackage{edale}
\usepackage[all]{xy}
\usepackage{mathrsfs}
\DeclareMathOperator{\Spec}{Spec}
\DeclareMathOperator{\Sp}{Sp}
\DeclareMathOperator{\Spa}{Spa}
\DeclareMathOperator{\Spf}{Spf}
\DeclareMathOperator{\cosq}{cosq}

\DeclareMathOperator{\Fr}{Fr}
\DeclareMathOperator{\Tr}{Tr}
\DeclareMathOperator{\Gal}{Gal}

\newcommand{\C}{\mathbb{C}}
\newcommand{\Z}{\mathbb{Z}}
\newcommand{\F}{\mathbb{F}}
\newcommand{\Q}{\mathbb{Q}}
\newcommand{\X}{\mathsf{X}}
\newcommand{\Y}{\mathsf{Y}}
\newcommand{\U}{\mathsf{U}}
\newcommand{\rig}{\mathrm{rig}}
\newcommand{\cf}{\textit{cf.\ }}
\newcommand{\word}[1]{{\bfseries\mathversion{bold}\relax #1}}

\makeatletter
\renewcommand{\maketitle}%
 {\@ifundefined{@title}{\relax}%
  {%
   \vspace*{30pt}\begin{center}\Large\word{\@title}\end{center}%
   \@ifundefined{@author}{\relax}{\begin{center}\large\@author\end{center}}%
   \thispagestyle{plain}%
  }%
 }
\makeatother

\title{On the action of the Weil group on the $\ell$-adic
 cohomology of rigid spaces over local fields}
\author{Yoichi Mieda}
\begin{document}
\maketitle

\begin{firstfootnote}
 Graduate School of Mathematical Sciences, the University of Tokyo,
 3--8--1 Komaba, Meguro-ku, Tokyo 153--8914, Japan.

 E-mail address: \texttt{edale@ms.u-tokyo.ac.jp}

\end{firstfootnote}

\begin{abstract}
 We investigate the action of the Weil group on the compactly supported $\ell$-adic cohomology groups
 of rigid spaces over local fields.
 We prove that every eigenvalue of the action is a Weil number when either a rigid space is smooth
 or the characteristic of the base field is equal to $0$.
 Since a smooth rigid space is locally isomorphic to the Raynaud generic fiber of an algebraizable
 formal scheme, we can use alterations and an analogue of the weight spectral sequence in the smooth case. 
 In the general case, we use the continuity theorem of Huber (\cite{Huber-JAG-!}, \cite{Huber-comparison}),
 which requires the restriction on the characteristic of the base field.
\end{abstract}

\section{Introduction}
Let $K$ be a complete discrete valuation field with finite residue field $\F_q$ and
$\overline{K}$ a separable closure of $K$.
We denote by $\Fr_q$ the geometric Frobenius element (the inverse of the $q$th power map)
in $\Gal(\overline{\F}_q/\F_q)$.
The Weil group $W_K$ of $K$ is defined as the inverse image of the subgroup
$\langle \Fr_q\rangle\subset\Gal(\overline{\F}_q/\F_q)$ by the canonical map
$\Gal(\overline{K}/K)\longrightarrow \Gal(\overline{\F}_q/\F_q)$.
For $\sigma\in W_K$, let $n(\sigma)$ be the integer such that the image of $\sigma$ in $\Gal(\overline{\F}_q/\F_q)$
is $\Fr_q^{n(\sigma)}$. Put $W_K^+=\{\sigma\in W_K\mid n(\sigma)\ge 0\}$.

Let $\X$ be a rigid space over $K$. 
In this paper, we consider the action of $W_K$ on the compactly supported $\ell$-adic cohomology group 
$H^i_c(\X\otimes_K\overline{K},\Q_\ell)$ (\cf \cite{Huber-book}, \cite{Huber-comparison}), 
where $\ell$ is a prime number which does not divide $q$. 
Our main theorem is the following:
 
\begin{thm}[][Theorem \ref{MainTheorem1}, Theorem \ref{MainTheorem2}]{Theorem}
 Let $\X$ be a quasi-compact separated rigid space over $K$. Assume one of the following conditions:
\begin{itemize}
 \item The rigid space $\X$ is smooth over $K$.
 \item The characteristic of $K$ is equal to $0$.
\end{itemize}
 Then for any $\sigma \in W_K^+$, every eigenvalue $\alpha\in \overline{\Q}_\ell$ of its action on
 $H^i_c(\X\otimes_K\overline{K},\Q_\ell)$
 is an algebraic integer. Moreover, there exists a non-negative integer $m$ such that for any isomorphism
 $\iota\colon \overline{\Q}_\ell\yrightarrow{\sim}\C$, the absolute value $\lvert \iota(\alpha)\rvert$ is
 equal to $q^{n(\sigma)\cdot m/2}$.
\end{thm}
Note that $H^i_c(\X\otimes_K\overline{K},\Q_\ell)$ is known to be a finite-dimensional $\Q_\ell$-vector space
when one of the above conditions is satisfied (\cite[Theorem 3.1]{Huber-comparison}). 

When $X$ is a scheme separated of finite type over $K$, the corresponding property is proven by Ochiai 
(\cite[Proposition 2.1]{Ochiai}).
Furthermore the $\ell$-independence of the alternating sum of the traces 
\[
 \sum_{i=0}^{2\dim X}(-1)^i\Tr\bigl(\sigma_*;H^i_c(X\otimes_K\overline{K},\Q_\ell)\bigr)
\]
is obtained there (\cite[Theorem 2.4]{Ochiai}). See also \cite{TSaito}, 
which treats the composite action of an element of $W_K$ and an algebraic correspondence.

We sketch the outline of the proof. When $\X$ is smooth over $K$, it is locally isomorphic to the
Raynaud generic fiber of an algebraizable formal scheme. In \S \ref{section-algebraization}, we briefly recall
this fact. By using the techniques of alterations and cohomological descent described in 
\S \ref{section-cohomological-descent}, the smooth case is reduced to the case where $\X$ is the Raynaud generic
fiber of the completion of a strictly semistable scheme over $\mathcal{O}_K$ along its special fiber.
In this case, we use an analogue of the weight spectral sequence to show the theorem. The smooth case is treated in \S
\ref{section-smoothcase}.
Assuming that the characteristic of $K$ is $0$, we prove the theorem for a general $\X$ by induction on $\dim \X$.
In this process, we need the continuity theorem of Huber (\cite{Huber-JAG-!}, \cite{Huber-comparison}).
The general case is treated in \S \ref{section-generalcase}.

\bigskip

\noindent\textbf{Acknowledgements}\ \,
The author would like to thank Tetsushi Ito for stimulating discussions, Takeshi Saito for pointing out
a mistake in the previous version of this paper. 
He is also grateful to his advisor Tomohide Terasoma.
He was supported by the Japan Society for the Promotion of Science Research Fellowships for
Young Scientists.

\bigskip

\noindent\textbf{Notation and Conventions}\ \,
Throughout this paper, $K$ is a complete discrete valuation field with finite residue field $\F_q$.
We denote the ring of integers of $K$ by $\mathcal{O}_K$.  
Fix a separable closure $\overline{K}$ of $K$ and write $X_{\overline{K}}$ for $X\otimes_K\overline{K}$.

We denote a scheme by an ordinary italic letter such as $X$,
a formal scheme by a calligraphic letter such as $\mathcal{X}$, and 
a rigid space by a sans serif letter such as $\X$.
For a scheme $X$ of finite type over $\mathcal{O}_K$,
$X_{s}$ (resp.\ $X_{\overline{s}}$) denotes its special fiber (resp.\ geometric special fiber) and
$X^{\wedge}$ the completion of $X$ along its special fiber.

For a formal scheme $\mathcal{X}$ over $\Spf \mathcal{O}_K$, we write $\mathcal{X}^\rig$ for its 
Raynaud generic fiber. It is the analytic adic space $d(\mathcal{X})$ in \cite[1.9]{Huber-book}.

We consider rigid spaces as adic spaces of finite type over $\Spa (K,\mathcal{O}_K)$.
We write $\Sp A$ for the adic space $\Spa (A,A^\circ)$, where $A$ is a topologically finitely generated $K$-algebra.

\section{Algebraization}\label{section-algebraization}
\subsection{}
In this section we see that a smooth rigid space over $K$ is locally isomorphic to the Raynaud generic fiber of 
an algebraizable formal scheme, i.e., the completion of a scheme of finite type over $\mathcal{O}_K$
along its special fiber. 

\subsection{}
Here we recall some definitions about topological rings introduced in \cite{Huber-valuation},
\cite{Huber-generalization}.
A topological ring is said to be {\slshape adic} if it has an ideal of definition, i.e., an ideal $I$ such that
$\{I^n\}_{n\ge 1}$ is a fundamental system of neighborhoods of $0$.
A topological ring $A$ is said to be {\slshape f-adic} if it has an open subring $A_0$ which has a finitely generated
ideal of definition. Such an open subring $A_0$ is called a {\slshape ring of definition}.
An f-adic ring is said to be {\slshape Tate} if it has a topologically nilpotent unit.

Let $A$ be a complete Tate ring. For a subset $U$ of $A$, we write $U\langle T_1,\ldots,T_n\rangle$ for the subset
\[
 \biggl\{\sum_{I\in (\Z_{\ge 0})^n} a_IT^I\in A[[T_1,\ldots,T_n]]\biggm| a_I\in U, \text{$(a_I)_{I\in (\Z_{\ge 0})^n}$ 
 converges to $0$ in $A$}\biggr\}
\]
of $A[[T_1,\ldots,T_n]]$. We equip $A\langle T_1,\ldots,T_n\rangle$ with the group topology such that
$U\langle T_1,\ldots,T_n\rangle$ for neighborhoods $U$ of $0\in A$ form a fundamental system
of neighborhoods of $0$. Then $A\langle T_1,\ldots,T_n\rangle$ is a complete Tate ring.
 It is easy to see that $A\langle T_1,\ldots,T_n\rangle\langle T_{n+1},\ldots,T_m\rangle$ is naturally isomorphic to
 $A\langle T_1,\ldots,T_m\rangle$ as a topological ring. 

\subsection{Example}
{\upshape A complete discrete valuation field $K$ with its natural topology is a complete Tate ring.
The ring of integers $\mathcal{O}_K$ is a ring of definition and a uniformizer of $K$ is
a topologically nilpotent unit of $K$.
The topological ring $K\langle T_1,\ldots,T_n\rangle$ is the same one as that in \cite{BGR}.}

\subsection{Lemma}\label{algebraizable}
Let $A$ be a complete Tate ring and $f_1,\ldots,f_n$ elements of $A\langle T_1,\ldots,T_n\rangle$.
Assume that the image of 
$\det (\partial f_i/\partial T_j)_{1\le i,j\le n}$ in $A\langle T_1,\ldots,T_n\rangle/(f_1,\ldots,f_n)$
is a unit. Then there exists an open neighborhood $V$ of
$0\in A\langle T_1,\ldots,T_n\rangle$ such that any $g_1,\ldots,g_n$ with $g_i\in f_i+V$ satisfy
the following two conditions:
\begin{itemize}
 \item The two topological rings $A\langle T_1,\ldots,T_n\rangle/(f_1,\ldots,f_n)$ and 
       $A\langle T_1,\ldots,T_n\rangle/(g_1,\ldots,g_n)$ are (topologically) $A$-isomorphic.
 \item The image of $\det (\partial g_i/\partial T_j)_{1\le i,j\le n}$ in 
       $A\langle T_1,\ldots,T_n\rangle/(g_1,\ldots,g_n)$ is a unit.
\end{itemize}

\begin{prf}
 This Lemma is due to \cite[ Proposition 1.7.1 ii) $\Longrightarrow$ iii)]{Huber-book}. 
 Note that the noetherian assumption \textit{loc.~cit.}\ is not necessary when $A$ is Tate.
 Indeed, in the proof given there, (3) is the only step where the noetherian assumption is used.
 In our case we may derive the bijectivity of 
 $f\colon B^{\triangleright}\longrightarrow B^{\triangleright}$ from that of
 $f\colon B_0\longrightarrow B_0$ and the equality $B^{\triangleright}=B_0[1/\pi]$, where $\pi$ is a topologically
 nilpotent unit of $A^{\triangleright}$.
 Nevertheless we use this lemma only for strongly noetherian rings $A$. 
\end{prf}

\subsection{Corollary}\label{Raynaud}
Let $\X$ be a smooth rigid space over $K$ and $x$ a point of $\X$.
Then there exist an open neighborhood $\U$ of $x$, 
a scheme $X$ of finite type over $\mathcal{O}_K$ with smooth generic fiber, and a $K$-isomorphism 
$\U\cong (X^{\wedge})^{\rig}$.

\begin{prf}
 Since $\X$ is smooth over $K$, there exists an open neighborhood $\U$ of $x$
 which is isomorphic to $\Sp K\langle T_1,\ldots,T_n\rangle/(f_1,\ldots,f_m)$ for
 some $m\le n$ and $f_1,\ldots,f_m\in K\langle T_1,\ldots,T_n\rangle$ such that the image of 
 $\det (\partial f_i/\partial T_j)_{1\le i,j\le m}$ in $K\langle T_1,\ldots,T_n\rangle/(f_1,\ldots,f_m)$ is a unit.
 Put $A=K\langle T_{m+1},\ldots,T_n\rangle$ and apply Lemma \ref{algebraizable}. Then we have an open neighborhood $V$
 of $0\in A\langle T_1,\ldots,T_m\rangle=K\langle T_1,\ldots,T_n\rangle$ satisfying the two conditions above.
 Since $K[T_1,\ldots,T_n]$ is dense in $K\langle T_1,\ldots,T_n\rangle$, there exist 
 $g_1,\ldots,g_m\in K[T_1,\ldots,T_n]$ such that the two topological rings
 $K\langle T_1,\ldots,T_n\rangle/(f_1,\ldots,f_m)$ and $K\langle T_1,\ldots,T_n\rangle/(g_1,\ldots,g_m)$ are
 $A$-isomorphic (\textit{a fortiori} $K$-isomorphic) and $\Delta=\det (\partial g_i/\partial T_j)_{1\le i,j\le m}$
 is invertible in $K\langle T_1,\ldots,T_n\rangle/(g_1,\ldots,g_m)$.
 Clearly we may assume that $g_i\in \mathcal{O}_K[T_1,\ldots,T_n]$.
 Let $\pi$ be a uniformizer of $K$.
 Since $x(\Delta)\neq 0$ for every $x\in \U=\Sp K\langle T_1,\ldots,T_n\rangle/(g_1,\ldots,g_m)$,
 there exists an integer $k$ such that $x(\Delta)\ge x(\pi^k)$ for every $x\in \U$ (\cite[Lemma 3.11]{Huber-valuation}).
 On the other hand, we have an isomorphism
 \[
 \bigl(\mathcal{O}_K[T_1,\ldots,T_{n+1}]/(g_1,\ldots,g_m,T_{n+1}\Delta-\pi^k)\bigr)^{\wedge}[1/\pi]
 \cong K\langle T_1,\ldots,T_{n+1}\rangle/(g_1,\ldots,g_m,T_{n+1}\Delta-\pi^k),
 \]
 where ${}^\wedge$ denotes the $\pi$-adic completion. 
 Thus we conclude that $\U\cong (X^{\wedge})^\rig$, where 
 $X=\Spec \mathcal{O}_K[T_1,\ldots,T_{n+1}]/(g_1,\ldots,g_m,T_{n+1}\Delta-\pi^k)$.
 Since 
 \[
  \Spec K[T_1,\ldots,T_{n+1}]/(g_1,\ldots,g_m,T_{n+1}\Delta-\pi^k)=\Spec K[T_1,\ldots,T_n,\Delta^{-1}]/(g_1,\ldots,g_m)
 \]
 is smooth over $K$, the corollary is proven.
\end{prf}

\subsection{}
When $\X\cong (X^{\wedge})^{\rig}$, the compactly supported $\ell$-adic cohomology groups of $\X$ can be
calculated by the nearby cycle of $X$:

\subsection{Theorem}\label{GAGA}
Let $X$ be a scheme which is separated of finite type over 
$\mathcal{O}_K$ and $\X$ the rigid space $(X^{\wedge})^{\rig}$.
Then $H^i_c(\X_{\overline{K}},\Q_\ell)$ is canonically isomorphic to $H_c^i(X_{\overline{s}},R\psi\Q_\ell)$. Moreover,
the isomorphism is $W_K$-equivariant.

\begin{prf}
 Apply \cite[Theorem 5.7.6]{Huber-book} to $X\otimes \mathcal{O}_{\overline{K}}$, where 
 $\mathcal{O}_{\overline{K}}$ is
 the integral closure of $\mathcal{O}_K$ in $\overline{K}$. The $W_K$-equivariantness is obvious.
\end{prf}

\section{Cohomological descent}\label{section-cohomological-descent}
\subsection{}
In the present section, we provide several results on the technique of cohomological descent
for $\ell$-adic cohomology of rigid spaces. Here we will use the terminology in \cite[Section 5]{Hodge3}
(see also \cite[Expos\'e V\textsuperscript{bis}]{SGA4-II}).

\subsection{Proposition}\label{algebraization-coh-descent}
Let $\X$ be a quasi-compact separated rigid space which is smooth over $K$.
Then we can construct an \'etale hypercovering $\U_\bullet\longrightarrow X$
such that $\U_n\cong (U_n^{\wedge})^\rig$ for some schemes $U_n$
which is separated of finite type over $\mathcal{O}_K$.
Moreover we have a following spectral sequence:
\[
 E_1^{-i,j}=H^j_c(\U_{i\overline{K}},\Q_\ell) \Longrightarrow H^{-i+j}_c(\X_{\overline{K}},\Q_\ell).
\]

\begin{prf}
 This is an easy consequence of Corollary \ref{Raynaud}.
 The existence of spectral sequence can be shown as in \cite[Remark 5.5.12 ii)]{Huber-book}.
\end{prf}

\subsection{Proposition}\label{cohomological-descent}
Let $\X$ and $\Y$ be rigid spaces and $f\colon \Y\longrightarrow \X$ a proper surjection.
Then it is a morphism universally of cohomological descent relative to torsion sheaves.

\begin{prf}
 This is due to the proper base change theorem for torsion sheaves, 
 which is proven in \cite[Theorem 4.4.1]{Huber-book} (\cf \cite[Expos\'e V\textsuperscript{bis}, Corollary 4.1.6]{SGA4-II} for the Betti cohomology and \textit{loc.~cit.}\ Proposition 4.3.2 for the \'etale cohomology of schemes).
\end{prf}

\subsection{Proposition}\label{alteration-coh-descent}
Let $X$ be a scheme of finite type over $\mathcal{O}_K$ and $n$ a positive integer. 
Then we can construct a system of field extensions $K\subset K_0\subset \cdots \subset K_n$ 
and an $n$-truncated proper hypercovering $Y_\bullet\longrightarrow X\otimes_{\mathcal{O}_K}\mathcal{O}_{K_n}$
satisfying the following conditions:
\begin{itemize}
 \item The scheme $Y_i$ is an $\mathcal{O}_{K_n}$-scheme and every structure morphism 
       is an $\mathcal{O}_{K_n}$-morphism.
 \item For every $i$ there exists a scheme $Y'_i$ which is strictly semistable (\cf \cite[section 1]{TSaito})
       over $\mathcal{O}_{K_i}$ such that $Y_i\cong Y'_i\otimes_{\mathcal{O}_{K_i}}\mathcal{O}_{K_n}$.
\end{itemize}
Moreover, the associated semi-simplicial rigid space $(Y_{\bullet}^{\wedge})^\rig$ is an $n$-truncated
proper hypercovering of $(X^{\wedge}\otimes_{\mathcal{O}_K}\mathcal{O}_{K_n})^\rig
=(X^{\wedge})^\rig\otimes_KK_n$. 

\begin{prf}
 All the assertions except the last are well-known consequences of \cite{deJong}. 
 The last one follows from Proposition \ref{cohomological-descent} and the two lemmas below.
\end{prf}

\subsection{Lemma}
The functor $X\longmapsto (X^{\wedge})^\rig$ from the category of schemes of finite type over $\mathcal{O}_K$
to the category of rigid spaces over $K$ preserves any finite projective limit (therefore the functor
commutes with $\cosq_n$).

\begin{prf}
 Since each of the categories has fiber products and the final object, it is sufficient to prove that the functor
 commutes with fiber products. It is known that the functor $^\rig$ commutes with fiber products 
 (\cite[Corollary 4.6]{BL1}). We will show that the functor $^{\wedge}$ commutes with fiber products.
 We may work locally: let $A$ be a finitely generated $\mathcal{O}_K$-algebra and $B$, $C$ be finitely
 generated $A$-algebras. 
 We should prove $(B\otimes_A C)^{\wedge}\cong B^{\wedge}\widehat{\otimes}_{A^\wedge}C^{\wedge}$, 
 which can be found in \cite[Chap.~0 (7.7.1)]{EGA1}.
\end{prf}

\subsection{Lemma}
Let $X$ and $Y$ be schemes of finite type over $\mathcal{O}_K$ and $f\colon Y\longrightarrow X$ a proper surjection.
Then the morphism $(f^{\wedge})^\rig\colon (Y^{\wedge})^\rig\longrightarrow (X^{\wedge})^\rig$ is also a
proper surjection.

\begin{prf}
 The lemma seems well-known, but we include its proof in the context of adic spaces. 
 First recall some general facts on adic spaces. Let $L$ be a non-archimedean field 
 (\cite[Definition 1.1.3]{Huber-book}) and $L^+$ be a valuation ring of $L$.
 For a scheme $X$ of finite type over $\mathcal{O}_K$, consider a $\Spa (K,K^+)$-morphism 
 $i\colon \Spa(L,L^+)\longrightarrow (X^{\wedge})^\rig$. Then we have the morphism of formal schemes
 $i_1\colon \Spf L^+\longrightarrow X^{\wedge}$ and that of schemes $i_2\colon \Spec L^+\longrightarrow X$.
 The maps $i\longmapsto i_1$ and $i_1\longmapsto i_2$ are clearly injective.
 Since every point of $(X^{\wedge})^\rig$ is analytic, the image of the generic point of $\Spec L^+$ under $i_2$
 lies on the generic fiber of $X$.
 Furthermore, for every $\mathcal{O}_K$-morphism $f\colon Y\longrightarrow X$, the following three sets are
 naturally identified:
 \begin{itemize}
  \item the set of morphisms of adic spaces $i'\colon \Spa(L,L^+)\longrightarrow (Y^{\wedge})^\rig$ such that
	$i=(f^{\wedge})^\rig\circ i'$;
  \item the set of morphisms of formal schemes $i'_1\colon \Spf L^+\longrightarrow Y^{\wedge}$ such that
	$i_1=f^{\wedge}\circ i'_1$;
  \item the set of morphisms of schemes $i'_2\colon \Spec L^+\longrightarrow Y$ such that
	$i_2=f\circ i'_2$.
 \end{itemize}
 Indeed, by local consideration, we can easily see that the two maps $i'\longmapsto i'_1$ and
 $i'_1\longmapsto i'_2$ are bijective. The adicness of the morphism
 $f^{\wedge}\colon Y^{\wedge}\longrightarrow X^{\wedge}$ is crucial for the second bijectivity.

 Now we will show the properness of $(f^{\wedge})^\rig$ by the valuative criterion (\cite[Corollary 1.3.9, Lemma 1.3.10]{Huber-book}). Let $L_1^+$ be another valuation ring of $L$. Suppose that a commutative diagram
 \[
  \xymatrix{
 \Spa(L,L_1^+)\ar[r]^-{i'}\ar[d]& (Y^{\wedge})^\rig\ar[d]^-{(f^{\wedge})^\rig}\\
 \Spa(L,L^+)\ar[r]^-{i}& (X^{\wedge})^\rig
 }
 \]
 is given. We should prove that there exists a unique morphism $i''\colon \Spa(L,L^+)\longrightarrow (Y^{\wedge})^\rig$
 that makes the following diagram commutative:
 \[
  \xymatrix{
 \Spa(L,L_1^+)\ar[r]^-{i'}\ar[d]& (Y^{\wedge})^\rig\ar[d]^-{(f^{\wedge})^\rig}\\
 \Spa(L,L^+)\ar@{-->}[ru]^-{i''}\ar[r]^-{i}& (X^{\wedge})^\rig\lefteqn{.}
 }
 \]
 By the consideration above, this is equivalent to proving
 that there exists a unique morphism $i''_2\colon \Spec L^+\longrightarrow Y$ that makes the following diagram
 commutative:
  \[
  \xymatrix{
 \Spec L_1^+\ar[r]^-{i'_2}\ar[d]& Y\ar[d]^-{f}\\
 \Spec L^+\ar@{-->}[ru]^-{i''_2}\ar[r]^-{i_2}& X\lefteqn{.}
 }
 \]
 This immediately follows from the valuative criterion for schemes.

 Next we will give a proof of the surjectivity of $(f^{\wedge})^\rig$. Let $x$ be a point of $(X^{\wedge})^\rig$
 and $i\colon \Spa (L,L^+)\longrightarrow (X^{\wedge})^\rig$ be a morphism which maps the closed point of
 $\Spa (L,L^+)$ to $x$. We may assume that the field $L$ is algebraically closed. We get the morphism of
 schemes $i_2\colon \Spec L^+\longrightarrow X$.
 Since $f$ is surjective and $L$ is algebraically closed, we can find a morphism $k\colon \Spec L\longrightarrow Y$
 such that the composite $f\circ k$ coincides with the restriction of $i_2$.
 By the properness of $f$, there exists a unique morphism $i_2'\colon \Spec L^+\longrightarrow Y$
 such that $i_2=f\circ i_2'$ and $i'_2\vert_{\Spec L}=k$.
 Let $i'\colon \Spa (L,L^+)\longrightarrow ({Y^{\wedge}})^\rig$ be 
 the corresponding morphism and $y$ be the image of the closed point of $\Spa (L,L^+)$ under $i'$.
 It is clear that the image of $y$ under $(f^{\wedge})^\rig$ coincides with $x$.
\end{prf}

\section{Smooth case}\label{section-smoothcase}
\subsection{}
In this section, we will give a proof of the main theorem for smooth rigid spaces.

\subsection{Theorem}\label{MainTheorem1}
Let $\X$ be a quasi-compact separated rigid space which is smooth over $K$.
Then for any $\sigma \in W_K^+$, every eigenvalue $\alpha\in \overline{\Q}_\ell$ of its action on
$H^i_c(\X_{\overline{K}},\Q_\ell)$
is an algebraic integer. Moreover, there exists a non-negative integer $m$
such that for any isomorphism $\iota\colon \overline{\Q}_\ell\yrightarrow{\sim}\C$, 
the absolute value $\lvert \iota(\alpha)\rvert$ is equal to $q^{n(\sigma)\cdot m/2}$.

\subsection{Lemma}\label{field-ext}
To prove Theorem \ref{MainTheorem1}, we have only to show that 
for some finite extension $L$ of $K$ the action of $W_L^+$ on 
$H^i_c(\X_{\overline{L}},\Q_\ell)$ satisfies the assertion in Theorem \ref{MainTheorem1}.

\begin{prf}
Since an inseparable extension affects neither the Weil group nor the \'etale cohomology, we may assume that
the extension $L/K$ is separable and that $L$ is a subfield of $\overline{K}$.
Write $e$ for the ramification index of $L/K$ and $f$ for the degree of the extension
of the residue field of $L/K$. Let $n'\colon W_L\longrightarrow \Z$ be the map $n$ for $L$ defined as in \S 1.
Then it is equal to the restriction of $n\colon W_K\longrightarrow \Z$ multiplied by $1/f$.
Take an element $\sigma$ of $W_K^+$ and an eigenvalue $\alpha$
of the action of $\sigma$ on $H^i_c(\X_{\overline{K}},\Q_\ell)$. Then we have $\sigma^{ef}\in W_L^+$.
By the assumption, the eigenvalue $\alpha^{ef}$ of the action of $\sigma^{ef}$ on $H^i_c(\X_{\overline{K}},\Q_\ell)$
is an algebraic integer and there exists a non-negative integer $m$
such that $\lvert \iota(\alpha^{ef})\rvert=(q^f)^{n'(\sigma^{ef})\cdot m/2}$ for
every isomorphism $\iota\colon \overline{\Q}_\ell\yrightarrow{\sim} \C$. 
Since $n'(\sigma^{ef})=n(\sigma^{ef})/f=e\cdot n(\sigma)$, the eigenvalue $\alpha$ is also an algebraic integer
and $\lvert \iota(\alpha)\rvert=q^{n(\sigma)\cdot m/2}$ for every $\iota$.
\end{prf}

\subsection{Lemma}\label{reduce-opencov}
To prove Theorem \ref{MainTheorem1}, we may assume that $\X\cong (X^{\wedge})^\rig$ for some scheme $X$
which is separated of finite type over $\mathcal{O}_K$.

\begin{prf}
Take an \'etale hypercovering $\U_\bullet\longrightarrow \X$ as in Proposition \ref{algebraization-coh-descent}.
Then we have the following spectral sequence:
\[
 E_1^{-i,j}=H^j_c(\U_{i\overline{K}},\Q_\ell) \Longrightarrow H^{-i+j}_c(\X_{\overline{K}},\Q_\ell).
\]
Every eigenvalue $\alpha$ of the action of $\sigma\in W_K^+$ on $H^{\nu}_c(\X_{\overline{K}}, \Q_\ell)$ occurs 
as an eigenvalue of the action of $\sigma$ on $H_c^j(\U_{i\overline{K}},\Q_\ell)$ for some non-negative integers 
$i$, $j$ satisfying $-i+j=\nu$. 
Thus we have to show Theorem \ref{MainTheorem1} only for $\U_i$, which is isomorphic to 
$(U_i^{\wedge})^\rig$ for some scheme $U_i$. This completes the proof.
\end{prf}

\subsection{Lemma}\label{reduce-alteration}
To prove Theorem \ref{MainTheorem1}, we may assume that $\X\cong (X^{\wedge})^\rig$ for some scheme $X$
which is strictly semistable over $\mathcal{O}_K$.

\begin{prf}
By Lemma \ref{reduce-opencov}, we may assume that $\X=(X^{\wedge})^\rig$ for some scheme $X$ which is separated
of finite type over $\mathcal{O}_K$.
It is sufficient to consider the action of $W_K$ on $H_c^\nu(\X_{\overline{K}},\Q_\ell)$ for $0\le \nu\le 2\dim \X$, 
since $H_c^\nu(\X_{\overline{K}},\Q_\ell)=0$ otherwise (\cite[Corollary 1.8.8 and Proposition 5.5.8]{Huber-book}).
Put $n=2\dim \X$ and take an $n$-truncated proper hypercovering
$Y_\bullet\longrightarrow X\otimes_{\mathcal{O}_K}\mathcal{O}_{K_n}$ as in Proposition \ref{alteration-coh-descent}.
Then we have the proper hypercovering $Y_\bullet=\cosq_n(Y_\bullet)$ and
the following $W_{K_n}$-equivariant spectral sequence:
\[
 E_1^{i,j}=H_c^j(\Y_{i\overline{K}},\Q_\ell)\Longrightarrow H_c^{i+j}(\X_{\overline{K}},\Q_\ell),
\]
where $\Y_i=(Y_i^{\wedge})^\rig$.
Every eigenvalue $\alpha$ of the action of $\sigma\in W_{K_n}^+$ on $H^\nu_c(\X_{\overline{K}}, \Q_\ell)$ occurs 
as an eigenvalue of the action of $\sigma$ on $H^j_c(\Y_{i\overline{K}},\Q_\ell)$ for some non-negative integers $i$,
$j$ satisfying $i+j=\nu$. 
Since $i=\nu-j\le \nu\le 2\dim \X$, there exist a scheme $Y_i'$ which is strictly semistable over $\mathcal{O}_{K_i}$
and an isomorphism $Y_i\cong Y_i'\otimes_{\mathcal{O}_{K_i}}\mathcal{O}_{K_n}$.
Then we have a $W_{K_n}$-equivariant isomorphism 
$H^j_c(\Y_{i\overline{K}},\Q_\ell)\cong H^j_c(\Y'_{i\overline{K}},\Q_\ell)$, where $\Y'_i=(Y_i'^{\wedge})^\rig$.
Since $Y_i'$ is strictly semistable over $\mathcal{O}_{K_i}$,
the eigenvalue $\alpha$ satisfies the property stated in Theorem \ref{MainTheorem1} by the assumption.
Combining this with Lemma \ref{field-ext}, we may conclude the present lemma.
\end{prf}

\subsection{}
By Theorem \ref{GAGA} and Lemma \ref{reduce-alteration},
we reduce Theorem \ref{MainTheorem1} to the following proposition:

\subsection{Proposition}\label{ssred}
Let $X$ be a strictly semistable scheme over $\mathcal{O}_K$. 
Then for any $\sigma \in W_K^+$, every eigenvalue $\alpha\in \overline{\Q}_\ell$ of its action on
$H^i_c(X_{\overline{s}},R\psi\Q_\ell)$ is an algebraic integer. 
Moreover, there exists a non-negative integer $m$ such that for any isomorphism
$\iota\colon \overline{\Q}_\ell\yrightarrow{\sim}\C$, the absolute value $\lvert \iota(\alpha)\rvert$ is
equal to $q^{n(\sigma)\cdot m/2}$.

\begin{prf}
 Denote the irreducible components of $X_s$ by $D_1,\ldots,D_m$. For a non-empty subset $I$ of $\{1,\ldots,m\}$
 and a non-negative integer $k$, we put $D_I=\bigcap_{i\in I}D_i$ and $D^{(k)}=\coprod_{\lvert I\rvert=k+1}D_I$.
 These are smooth over $\F_q$. We use the spectral sequence
 \[
  E_1^{i,j}=\bigoplus_{k\ge\max (0,-i)}H_c^{j-2k}\bigl(D^{(i+2k)}_{\overline{s}},\Q_\ell(-k)\bigr)
 \Longrightarrow H_c^{i+j}(X_{\overline{s}},R\psi\Q_\ell).
 \]
 This is the spectral sequence associated with the monodromy filtration of $R\psi\Q_\ell$ (\cite{TSaito})
 and is called the weight spectral sequence when $X$ is proper over $\mathcal{O}_K$ (\cite{RaZi}).

 For every $\sigma\in W_K^+$, we have the following morphism of spectral sequences 
 (\cf \cite[proof of Lemma 3.2]{TSaito}):
 \[
  \xymatrix{
 E_1^{i,j}=\bigoplus_{k\ge\max (0,-i)}H_c^{j-2k}\bigl(D^{(i+2k)}_{\overline{s}},\Q_\ell(-k)\bigr)
 \ar@{=>}[r]\ar[d]^-{\bigl(\Fr_q^{n(\sigma)}\bigr)_*}& H_c^{i+j}(X_{\overline{s}},R\psi\Q_\ell)\ar[d]^-{\sigma_*}\\
 E_1^{i,j}=\bigoplus_{k\ge\max (0,-i)}H_c^{j-2k}\bigl(D^{(i+2k)}_{\overline{s}},\Q_\ell(-k)\bigr)
 \ar@{=>}[r]& H_c^{i+j}(X_{\overline{s}},R\psi\Q_\ell)\lefteqn{.}
 }
 \]
 Therefore every eigenvalue $\alpha$ of the action of $\Fr_q^{n(\sigma)}$ on $H^\nu_c(X_{\overline{s}},R\psi\Q_\ell)$
 occurs as an eigenvalue of the action of $\Fr_q^{n(\sigma)}$ on 
 $H_c^{j-2k}(D^{(i+2k)}_{\overline{s}},\Q_\ell(-k))$ for some $i$, $j$, $k$ satisfying
 $i+j=\nu$, $k\geq \max (0,-i)$.
 On the other hand, by the purity theorem of Deligne (\cite[Corollaire 3.3.4]{Weil2}) and 
 \textit{loc.~cit.}\ Corollaire 3.3.3, every eigenvalue $\alpha$ of the action of $\Fr_q^{n(\sigma)}$ on
 $H_c^{j-2k}(D^{(i+2k)}_{\overline{s}},\Q_\ell(-k))$ is an algebraic integer
 and there exists an integer $m$ with $2k\le m$ such that
 $\lvert\iota(\alpha)\rvert=q^{n(\sigma)\cdot m/2}$ for every isomorphism
 $\iota\colon \overline{\Q}_\ell\yrightarrow{\sim}\C$. Since $k\ge 0$, the integer $m$ is non-negative.
 This completes the proof.
\end{prf}

\section{General case}\label{section-generalcase}
\subsection{}
In this section, we assume that the characteristic of $K$ is equal to $0$.
Recall the following theorems by Huber, which are used to prove the main theorem.

\subsection{Theorem}\label{finitude}
Assume that the characteristic of $K$ is equal to $0$.
Let $\X$ be a quasi-compact separated rigid space over $K$.
Then the cohomology group $H^i_c(\X_{\overline{K}},\Q_\ell)$ is a finite-dimensional $\Q_\ell$-vector space.

\begin{prf}
 \cite[Theorem 3.1]{Huber-comparison} (see also \cite[Corollary 2.3]{Huber-JAG-!}).
\end{prf}

\subsection{Theorem}\label{continuity}
Let $\X$ be a separated rigid space over $K$ which is not necessarily quasi-compact, 
$\{\U_\lambda\}$ an open covering of $\X$ consisting of quasi-compact open subspaces
such that for every $\lambda_1$, $\lambda_2$ there exists $\lambda$ satisfying
$\U_{\lambda_1}\cup \U_{\lambda_2}\subset \U_{\lambda}$.
Then the canonical homomorphism
 \[
  \varinjlim_{\lambda}H^i_c(\U_{\lambda\overline{K}},\Q_\ell)\longrightarrow H^i_c(\X_{\overline{K}},\Q_\ell)
 \]
is an isomorphism.
 
\begin{prf}
 \cite[Proposition 2.1 (iv)]{Huber-comparison}. Note that this theorem is also valid when the characteristic of $K$
 is positive.
\end{prf}

\subsection{}
Now we can give a proof of our main theorem.

\subsection{Theorem}\label{MainTheorem2}
Assume that the characteristic of $K$ is equal to $0$.
Let $\X$ be a quasi-compact separated rigid space over $K$. 
Then for any $\sigma \in W_K^+$, every eigenvalue $\alpha\in \overline{\Q}_\ell$ of its action on
$H^i_c(\X_{\overline{K}},\Q_\ell)$
is an algebraic integer. Moreover, there exists a non-negative integer $m$ such that for any isomorphism
$\iota\colon \overline{\Q}_\ell\yrightarrow{\sim}\C$, the absolute value $\lvert \iota(\alpha)\rvert$ is
equal to $q^{n(\sigma)\cdot m/2}$.

\begin{prf}
 We proceed by induction on $\dim \X$. We may assume that $\X$ is reduced.
 Let $\mathsf{Z}$ be the singular locus of $\X$. It is a closed analytic subspace whose dimension is strictly
 less than $\dim \X$.
 Thus we have only to show our claim on $H^i_c(\U_{\overline{K}},\Q_\ell)$, where $\U=\X\setminus \mathsf{Z}$.
 Let $\{\U_\lambda\}$ be the set of quasi-compact open subspaces of $\U$ (note that $\U$ itself is rarely quasi-compact).
 Then we have the isomorphism 
 \[
  \varinjlim_{\lambda}H^i_c(\U_{\lambda\overline{K}},\Q_\ell)\yrightarrow {\sim} H^i_c(\U_{\overline{K}},\Q_\ell)
 \]
 by Theorem \ref{continuity}.
 Since every eigenvalue of the action of $\sigma\in W_K^+$ on $H^i_c(\U_{\overline{K}},\Q_\ell)$ occurs as
 an eigenvalue of the action of $\sigma$ on $H^i_c(\U_{\lambda\overline{K}},\Q_\ell)$ for some $\lambda$,
 and each $\U_\lambda$ is quasi-compact and smooth over $K$,
 the theorem follows from Theorem \ref{MainTheorem1}.
\end{prf}

\subsection{}
If we can prove Theorem \ref{finitude} for $K$ of positive characteristic, 
we can remove the restriction of the characteristic.


\begin{thebibliography}{99}
 \bibitem[BGR]{BGR} S.~Bosch, U.~Guntzer, R.~Remmert, \textit{Non-Archimedean analysis},
	      Grundlehren der Mathematischen Wissenschaften, 261. Springer-Verlag, Berlin, 1984.
 \bibitem[BL]{BL1} S.~Bosch, W.~L\"utkebohmert, \textit{Formal and rigid geometry I. Rigid spaces.}, Math.\ Ann.\ 295 
	     (1993), no.~2, 291--317. 
 \bibitem[dJ]{deJong} A.~J.~de Jong, \textit{Smoothness, semi-stability and alterations},
	     Inst.\ Hautes \'Etudes Sci.\ Publ.\ Math., No.~83 (1996), 51--93. 
 \bibitem[De1]{Weil2} P.~Deligne, \textit{La conjecture de Weil II}, Inst.\ Hautes \'Etudes Sci.\ Publ.\ Math.\ 
	      No.~52 (1980), 137--252.
 \bibitem[De2]{Hodge3} P.~Deligne, \textit{Theorie de Hodge III}, Inst.\ Hautes \'Etudes Sci.\ Publ.\ Math., No.~44
	      (1974), 5--77.
 \bibitem[Hu1]{Huber-valuation} R.~Huber, \textit{Continuous valuations}, Math.\ Z.\ 212 (1993), no.~3, 455--477.
 \bibitem[Hu2]{Huber-generalization} R.~Huber, \textit{A generalization of formal schemes and rigid analytic varieties},
	      Math.\ Z.\ 217 (1994), no.~4, 513--551.
 \bibitem[Hu3]{Huber-book} R.~Huber, \textit{\'Etale cohomology of rigid analytic varieties and adic spaces},
	      Aspects of Mathematics, E30. Friedr.\ Vieweg \& Sohn, Braunschweig, 1996.
 \bibitem[Hu4]{Huber-JAG-!} R.~Huber, \textit{A finiteness result for the compactly supported cohomology of rigid analytic varieties}, J.\ Algebraic Geom.\ 7 (1998), no.~2, 313--357.
 \bibitem[Hu5]{Huber-comparison} R.~Huber, \textit{A comparison theorem for $\ell$-adic cohomology}, 
	      Compositio Math.\ 112 (1998), no.~2, 217--235.
 \bibitem[O]{Ochiai}  T.~Ochiai, \textit{$\ell$-independence of the trace of monodromy}, Math.\ Ann.\ 315 (1999),
	    no.~2, 321--340.
 \bibitem[Sa]{TSaito} T.~Saito, \textit{Weight spectral sequences and independence of $\ell$}, J.\ Inst.\ Math.\ Jussieu 2 (2003), no.~4, 583--634.
 \bibitem[RZ]{RaZi} M.~Rapoport, Th.~Zink, \textit{\"Uber die lokale Zetafunktion von Shimuravariet\"aten. 
	     Monodromiefiltration und verschwindende Zyklen in ungleicher Charakteristik}, 
	     Invent.\ Math.\ 68 (1982), no.~1, 21--101.
 \bibitem[EGA1]{EGA1} A.~Grothendieck, \textit{Elements de geometrie algebrique I. Le langage des sch\'emas.},
	       Inst.\ Hautes \'Etudes Sci.\ Publ.\ Math.\ No.~4 (1960).
 \bibitem[SGA4-II]{SGA4-II} \textit{Th\'eorie des topos et cohomologie \'etale des sch\'emas}, Tome 2.
		 S\'eminaire de G\'eom\'etrie Alg\'ebrique du Bois-Marie 1963--1964 (SGA 4). 
		 Dirig\'e par M.~Artin, A.~Grothendieck et J.~L.~Verdier. 
		 Avec la collaboration de N.~Bourbaki, P.~Deligne et B.~Saint-Donat. 
		 Lecture Notes in Mathematics, Vol.~270. Springer-Verlag, Berlin-New York, 1972.
\end{thebibliography}
\end{document}